# Stabilization of 2D discrete switched systems with state delays under asynchronous switching

**Shipei Huang,** **Zhengrong Xiang***

School of Automation, Nanjing University of Science and Technology

Nanjing, 210094, People's Republic of China

*corresponding author, e-mail: xiangzr@mail.njust.edu.cn

**Abstract:** This paper is concerned with the problem of robust stabilization for a class of uncertain 2D discrete switched systems with state delays represented by a model of Roesser type, where the switching instants of the controller experience delays with respect to those of the system, and the parameter uncertainties are assumed to be norm-bounded. A state feedback controller is proposed to guarantee exponential stability for such 2D discrete switched systems, and the dwell time approach is utilized for the stability analysis and controller design. A numerical example is given to illustrate the effectiveness of the proposed method.

## 1. Introduction

Two-dimensional (2D) system consists of two variables in horizontal and vertical, which is different from 1D system. Such systems have received considerable attention over the past few decades due to their wide applications in many areas, such as multi-dimensional digital filtering, linear image processing, signal processing, and process control [1-3]. It is well known that 2D systems can be represented by different models such as the Roesser model, Fornasini–Marchesini model and Attasi model. The issues of stability analysis and control synthesis of 2-D discrete systems have been studied in [4-13].

It is noticed that the time-delay phenomenon exists widely in engineering and social systems, time delays frequently occur in practical systems and time delays are often the source of instability. There are many examples containing inherent delays in practical 2-D discrete systems. The stability of 2-D discrete systems with state delays has also been studied in [8-12].

On the other hand, switched systems have also been attracted considerable attention during the past several decades [14-26]. A switched system is a hybrid system which consists of a finite number of continuous-time or discrete-time subsystems and a switching signal specifying the switch between these subsystems. This class of systems has numerous applications in many fields, such as mechanical systems, the automotive industry, aircraft and air traffic control, switched power converters. So far, there are three methods many authors have been used to study such systems. The first method is common Lyapunov quadratic functions and the second is multiple Lyapunov functions. The last approach is different from them, that is dwell time approach. Recently, the dwell time approach is applied widely to deal with switched systems, see, for example, [22-26] and references therein.

However, switch phenomenon may also occur in practical 2D discrete systems, and studying of these 2D discrete switched systems will also be significant. There are a few reports on 2D discrete switched systems at present. Benzaouia, Hmamed, and Tadeo [27] firstly consider 2D switched systems with arbitrary switching sequences, and the process of switching is considered as a Markovian jumping one. In addition, they study the stabilizability problem of 2D discrete switched systems in [28], by using the common Lyapunov function approach and multiple Lyapunov function approach, two different sufficient conditions are proposed for the 2D discrete switched system such that the correspongding closed-loop system is asymptotic stability, respectively.

It is noticed that there inevitably exists asynchronous switching between the controller and the system in actual operation. Investigating the problems of stability and control synthesis for switched systems under asynchronous switching could be very significant. Many results on the issues of stability and control synthesis for the 1D switched systems under asynchronous switching have been obtained in [23-26], most of which are based on the dwell time approach. However, to the best of our knowledge, the problems of stability and stabilization for 2D switched systems with state delays under asynchronous switching via the dwell time approach have not been investigated to date, especially for the exponential stability problem of these systems, which motivates our present study.

In this paper, we are interested in designing a stabilizing controller for 2D discrete switched systems with state delays represented by a model of Roesser type under asynchronous switching such that the corresponding closed-loop system is exponentially stable. The dwell time approach is utilized for the stability analysis and controller design.

This paper is organized as follows. In Section 2, problem formulation and some necessary lemmas are given. In Section 3, based on the dwell time approach, stability and stabilization for 2D discrete switched systems with state delays are addressed, and sufficient conditions for the existence of a stabilizing controller for such 2D discrete switched systems under asynchronous switching are derived in terms of a set of matrix inequalities. A numerical example is provided to illustrate the effectiveness of the proposed approach in Section 4. Concluding remarks are given in Section 5.

**Notations**: Throughout this paper, the superscript "$T$" denotes the transpose, and the notation $X \geq Y \left( X > Y \right)$ means that matrix $X - Y$ is positive semi-definite (positive definite, respectively). $\|*\|$ denotes the Euclidean norm. $I$ represents identity matrix with appropriate dimension. $I_h$ is the identity matrix with $n_1$ appropriate dimension and $I_v$ is the identity matrix with $n_2$ appropriate dimension. $diag\{a_i\}$ denotes diagonal matrix with the diagonal elements $a_i, i = 1, 2, ..., n$. $X^{-1}$ denotes the inverse of $X$. The asterisk $*$ in a matrix is used to denote term that is induced by symmetry. The set of all nonnegative integers is represented by $Z_+$.

## 2. Problem formulation and preliminaries

Consider the following uncertain 2D discrete linear switched systems with state delays described by Roesser model:

$$\begin{bmatrix} x^h\left(i+1, j\right) \\ x^v\left(i, j+1\right) \end{bmatrix} = \hat{A}^{\sigma(m)} \begin{bmatrix} x^h\left(i, j\right) \\ x^v\left(i, j\right) \end{bmatrix} + \hat{A}_d^{\sigma(m)} \begin{bmatrix} x^h\left(i-d_h, j\right) \\ x^v\left(i, j-d_v\right) \end{bmatrix} + \hat{B}^{\sigma(m)} u\left(i, j\right) \tag{1}$$

where $x^h\left(i, j\right)$ is the horizontal state in $R^{n_1}$, $x^v\left(i, j\right)$ is the vertical state in $R^{n_2}$, $x\left(i, j\right)$ is the whole state in $R^n$ with $n = n_1 + n_2$. $u\left(i, j\right)$ is the control input. $m$ is the sum of $i$ and $j$. $\sigma\left(m\right)$ is a switching rule which takes its values in the finite set $\underline{N} := \{1, \cdots, N\}$, $N$ is the number of subsystems, and $\sigma\left(m\right) = k$ means that the $k-th$ subsystem is activated. $i$ and $j$ are integers in $Z_+$. $d_h$ and $d_v$ are constant delays along horizontal and vertical directions, respectively. $\hat{A}^{\sigma(m)}$, $\hat{A}_d^{\sigma(m)}$ are uncertain real-valued matrices with appropriate dimensions which are assumed to be of the form

$$\hat{A}^{\sigma(m)} = A^{\sigma(m)} + H^{\sigma(m)} F^{\sigma(m)}(i,j) E_1^{\sigma(m)}$$

$$\hat{A}_d^{\sigma(m)} = A_d^{\sigma(m)} + H^{\sigma(m)} F^{\sigma(m)}(i,j) E_2^{\sigma(m)}$$

$$\hat{B}^{\sigma(m)} = B^{\sigma(m)} + H^{\sigma(m)} F^{\sigma(m)}(i,j) E_3^{\sigma(m)} \tag{2}$$

with

$$A^{\sigma(m)} = \begin{bmatrix} A_{11}^{\sigma(m)} & A_{12}^{\sigma(m)} \\ A_{21}^{\sigma(m)} & A_{22}^{\sigma(m)} \end{bmatrix},\ A_d^{\sigma(m)} = \begin{bmatrix} A_{d11}^{\sigma(m)} & A_{d12}^{\sigma(m)} \\ A_{d21}^{\sigma(m)} & A_{d22}^{\sigma(m)} \end{bmatrix},$$

$$H^{\sigma(m)} = \begin{bmatrix} H_1^{\sigma(m)} \\ H_2^{\sigma(m)} \end{bmatrix},\ E_1^{\sigma(m)} = \begin{bmatrix} E_{11}^{\sigma(m)} \\ E_{12}^{\sigma(m)} \end{bmatrix},\ E_2^{\sigma(m)} = \begin{bmatrix} E_{21}^{\sigma(m)} \\ E_{22}^{\sigma(m)} \end{bmatrix},\ E_3^{\sigma(m)} = \begin{bmatrix} E_{31}^{\sigma(m)} \\ E_{32}^{\sigma(m)} \end{bmatrix}.$$

where matrices $A_{11}^{\sigma(m)} \in R^{n_1 \times n_1}$, $A_{12}^{\sigma(m)} \in R^{n_1 \times n_2}$, $A_{21}^{\sigma(m)} \in R^{n_2 \times n_1}$, $A_{22}^{\sigma(m)} \in R^{n_2 \times n_2}$, $A_{d11}^{\sigma(m)} \in R^{n_1 \times n_1}$, $A_{d12}^{\sigma(m)} \in R^{n_1 \times n_2}$, $A_{d21}^{\sigma(m)} \in R^{n_2 \times n_1}$, $A_{d22}^{\sigma(m)} \in R^{n_2 \times n_2}$, $H_1^{\sigma(m)}$, $H_2^{\sigma(m)}$, $E_{11}^{\sigma(m)}$, $E_{12}^{\sigma(m)}$, $E_{21}^{\sigma(m)}$, $E_{22}^{\sigma(m)}$, $E_{31}^{\sigma(m)}$, $E_{32}^{\sigma(m)}$ are constant matrices. $F^{\sigma(m)}(i,j)$ is an unknown matrix representing parameter uncertainty which satisfies

$$F^{\sigma(m)T}(i,j) F^{\sigma(m)}(i,j) \le I \tag{3}$$

The boundary conditions are defined by

$$\begin{aligned} &x^h(i,j) = h_{ij}, \quad \forall 0 \le j \le z_1,\ -d_h \le i \le 0 \\ &x^h(i,j) = 0, \quad \forall j > z_1, \quad -d_h \le i \le 0 \\ &x^v(i,j) = v_{ij}, \quad \forall 0 \le i \le z_2,\ -d_v \le j \le 0 \\ &x^v(i,j) = 0, \quad \forall i > z_2, \quad -d_v \le j \le 0 \end{aligned} \tag{4}$$

where $z_1 < \infty$ and $z_2 < \infty$ are positive integers, $h_{ij}$ and $v_{ij}$ are given vectors.

**Remark 1** In the paper [28], it has stated that since the 2D system causality imposes an increment depending only on $d = i + j$, the switch can be assumed to occur at $m = i + j$. In the paper, we can make the same assumption, and the switching sequence can be described as

$$(m_0, \sigma(m_0)), (m_1, \sigma(m_1)), \ldots, (m_\kappa, \sigma(m_\kappa)) \ldots \tag{5}$$

with $\kappa \in Z_+$, $m_\kappa = i_\kappa + j_\kappa$, $m_\kappa$ denotes the $\kappa - th$ switching instant and $\sigma(m_\kappa) \in \underline{N}$ denotes the value of switching signal.

However, in actual operation, there inevitably exists asynchronous switching between the controller and the system. Without loss of generality, we only consider the case that the switching instants of the controller experience delays with respect to the switching instants of the system. Let $\sigma'(m)$ denote the switching signal of the controller, the switching instants of the controller can be described as

$$m_1 + \Delta m_1,\ m_2 + \Delta m_2, \ldots,\ m_\kappa + \Delta m_\kappa, \ldots$$

where $\Delta m_\kappa < \inf(m_{\kappa+1} - m_\kappa)$, $\Delta m_\kappa$ represents the delayed period, and it is said to be mismatched period.

**Remark 2** Mismatched period $\Delta m_\kappa < \inf(m_{\kappa+1} - m_\kappa)$ guarantees that there always exists a period that the controller and the system operate synchronously, and this period is said to be matched period in the later section.

**Definition 1** System (1) is said to be exponentially stable under $\sigma(m)$ if for a given $z \ge 0$, there exists a positive constant $c$, such that

$$\sum_{i+j=D} \left\| x(i,j) \right\|^2 \le \eta e^{-c(D-z)} \sum_{i+j=z} \left\| x(i,j) \right\|_C^2 \tag{6}$$

holds for all $D \ge z$ and a positive constant $\eta$ and

$$\sum_{i+j=z} \left\| x(i,j) \right\|_C^2 \triangleq \sup_{\substack{-d_h < \theta_h \le 0, \\ -d_v < \theta_v \le 0}} \sum_{i+j=z} \left( \left\| x^h(i-\theta_h, j) \right\|^2 + \left\| x^v(i, j-\theta_v) \right\|^2 \right)$$

**Remark 3** From the definition 1, it is easy to see that when $z$ is given, $\sum_{i+j=z} \left\| x(i,j) \right\|_C^2$ will be bounded, and $\sum_{i+j=D} \left\| x(i,j) \right\|^2$ will tend to zero exponentially as $D$ goes to infinity, it implies that $\left\| x(i,j) \right\|$ tends to zero.

**Definition 2** Let $m_\kappa = i_\kappa + j_\kappa$ denotes the $\kappa$-th switching instant and $m_{\kappa+1} = i_{\kappa+1} + j_{\kappa+1}$ denotes the $\kappa+1$-th switching instant. Denote $\tau=\inf\left(m_{\kappa+1}-m_{\kappa}\right)$, $\tau$ is called the dwell time.

**Definition 3** For any $i+j=D\geq z=i_z+j_z$, let $N_{\sigma(m)}\left(z,D\right)$ denote the switching number of $\sigma\left(m\right)$ on an interval $\left(z,D\right)$. If

$$N_{\sigma(m)}\left(z,D\right)\leq N_0+\frac{D-z}{\tau_a} \tag{7}$$

holds for given $N_0\geq 0$, $\tau_a\geq 0$, then the constant $\tau_a$ is called the average dwell time and $N_0$ is the chatter bound.

**Lemma 1** [29] For a given matrix $S=\begin{bmatrix} S_{11} & S_{12} \\ S_{12}^T & S_{22} \end{bmatrix}$, where $S_{11}$, $S_{22}$ are square matrices, then the following conditions are equivalent.

**(i)** $S<0$;

**(ii)** $S_{11}<0, S_{22}-S_{12}^T S_{11}^{-1} S_{12}<0$;

**(iii)** $S_{22}<0$, $S_{11}-S_{12}S_{22}^{-1}S_{12}^T<0$;

**Lemma 2** [30] Let $U$, $V$, $W$ and $X$ be real matrices of appropriate dimensions with $X$ satisfying $X=X^T$, then for all $V^TV\leq I$, $X+UVW+W^TV^TU^T<0$, if and only if there exists a scalar $\varepsilon$ such that $X+\varepsilon UU^T+\varepsilon^{-1}W^TW<0$.

# 3. Main results

### 3.1. Stability analysis

In this subsection, we first focus on the problem of stability analysis for the 2D discrete systems with state delays.

**Lemma 3** Consider the following 2D discrete system with state delays

$$\begin{bmatrix} x^h\left(i+1,j\right) \\ x^v\left(i,j+1\right) \end{bmatrix}=A\begin{bmatrix} x^h\left(i,j\right) \\ x^v\left(i,j\right) \end{bmatrix}+A_d\begin{bmatrix} x^h\left(i-d_h,j\right) \\ x^v\left(i,j-d_v\right) \end{bmatrix} \tag{8}$$

where $A$ and $A_d$ are constant matrices with appropriate dimensions. The boundary conditions are as (4). For a given positive constant $\alpha < 1$, if there exist positive definite symmetric matrices $P = diag\{P_h, P_v\}$, $Q = diag\{Q_h, Q_v\}$ with appropriate dimensions, such that

$$\begin{bmatrix} Q-\alpha P & 0 & A^T P \\ * & -\Lambda_1 Q & A_d^T P \\ * & * & -P \end{bmatrix} < 0 \tag{9}$$

where $\Lambda_1 = diag\left\{\alpha^{d_h} I_h, \alpha^{d_v} I_v\right\}$.

then, along the trajectory of systems (8), there holds the following inequality

$$\sum_{i+j=D} V(i,j) < \alpha^{D-D'} \sum_{i+j=D'} V(i,j) \tag{10}$$

where $D' \geq z$.

**Proof** Consider Lyapunov-Krasovskii functional candidate

$$V(x(i,j)) = V^h\left(x^h(i,j)\right) + V^v\left(x^v(i,j)\right) \tag{11}$$

where

$$V^h\left(x^h(i,j)\right) = \sum_{g=1}^{2} V_g^h\left(x^h(i,j)\right)$$

$$V_1^h\left(x^h(i,j)\right) = x^h(i,j)^T P_h x^h(i,j)$$

$$V_2^h\left(x^h(i,j)\right) = \sum_{r=i-d_h}^{i-1} x^h(r,j)^T Q_h x^h(r,j)\alpha^{i-r-1}$$

$$V^v\left(x^v(i,j)\right) = \sum_{g=1}^{2} V_g^v\left(x^v(i,j)\right)$$

$$V_1^v\left(x^v(i,j)\right) = x^v(i,j)^T P_v x^v(i,j)$$

$$V_2^v\left(x^v(i,j)\right) = \sum_{t=j-d_v}^{j-1} x^v(i,s)^T Q_v x^v(i,s)\alpha^{j-t-1}$$

Then we have

$$V^h\left(x^h(i+1,j)\right) - \alpha V^h\left(x^h(i,j)\right) + V^v\left(x^v(i,j+1)\right) - \alpha V^v\left(x^v(i,j)\right)$$

$$= \sum_{g=1}^{2}\left[V_g^h\left(x^h\left(i+1,j\right)\right)-\alpha V_g^h\left(x^h\left(i,j\right)\right)\right]+\sum_{g=1}^{2}\left[V_g^v\left(x^v\left(i,j+1\right)\right)-\alpha V_g^v\left(x^v\left(i,j\right)\right)\right] \quad (12)$$

with

$$V_1^h\left(x^h\left(i+1,j\right)\right)-\alpha V_1^h\left(x^h\left(i,j\right)\right)$$

$$= x^h\left(i+1,j\right)^T P_h x^h\left(i+1,j\right)-\alpha x^h\left(i,j\right)^T P_h x^h\left(i,j\right)$$

$$V_2^h\left(x^h\left(i+1,j\right)\right)-\alpha V_2^h\left(x^h\left(i,j\right)\right)$$

$$= \sum_{r=i+1-d_h}^{i} x^h\left(r,j\right)^T Q_h x^h\left(r,j\right)\alpha^{i-r} - \alpha\sum_{r=i-d_h}^{i-1} x^h\left(r,j\right)^T Q_h x^h\left(r,j\right)\alpha^{i-r-1}$$

$$= x^h\left(i,j\right)^T Q_h x^h\left(i,j\right)-\alpha^{d_h} x^h\left(i-d_h,j\right)^T Q_h x^h\left(i-d_h,j\right)$$

$$V_1^v\left(x^v\left(i,j+1\right)\right)-\alpha V_1^h\left(x^v\left(i,j\right)\right)$$

$$= x^v\left(i,j+1\right)^T P_v x^v\left(i,j+1\right)-\alpha x^v\left(i,j\right)^T P_v x^v\left(i,j\right)$$

$$V_2^v\left(x^v\left(i,j+1\right)\right)-\alpha V_2^v\left(x^v\left(i,j\right)\right)$$

$$= \sum_{t=j+1-d_v}^{j} x^v\left(i,t\right)^T Q_v x^v\left(i,t\right)\alpha^{j-t} - \alpha\sum_{t=j-d_v}^{j-1} x^v\left(i,t\right)^T Q_v x^v\left(i,t\right)\alpha^{j-t-1}$$

$$= x^v\left(i,j\right)^T Q_v x^v\left(i,j\right)-\alpha^{d_v} x^v\left(i,j-d_v\right)^T Q_v x^v\left(i,j-d_v\right)$$

we can obtain the following equality

$$V^h\left(x^h\left(i+1,j\right)\right)-\alpha V^h\left(x^h\left(i,j\right)\right)+V^v\left(x^v\left(i,j+1\right)\right)-\alpha V^v\left(x^v\left(i,j\right)\right)$$

$$= \begin{bmatrix} \begin{bmatrix} x^h\left(i,j\right) \\ x^v\left(i,j\right) \end{bmatrix} \\ \begin{bmatrix} x^h\left(i-d_h,j\right) \\ x^h\left(i,j-d_v\right) \end{bmatrix} \end{bmatrix}^T \begin{bmatrix} \Phi_{11} & \Phi_{12} \\ \Phi_{12}^T & \Phi_{22} \end{bmatrix} \begin{bmatrix} \begin{bmatrix} x^h\left(i,j\right) \\ x^v\left(i,j\right) \end{bmatrix} \\ \begin{bmatrix} x^h\left(i-d_h,j\right) \\ x^h\left(i,j-d_v\right) \end{bmatrix} \end{bmatrix}$$

where

$$\Phi_{11} = Q-\alpha P + A^T P A,\ \Phi_{12} = A^T P A_d,\ \Phi_{22} = A_d^T P A_d - \Lambda_1 Q,\ \Lambda_1 = diag\left\{\alpha^{d_h} I_h, \alpha^{d_v} I_v\right\}$$

Applying lemma 1, (9) is equivalent to the following inequality

$$\begin{bmatrix} \Phi_{11} & \Phi_{12} \\ \Phi_{12}^T & \Phi_{22} \end{bmatrix} < 0 \quad (13)$$

For simplicity, we define

$$V^h(i,j)=V^h\left(x^h(i,j)\right),\ V^v(i,j)=V^v\left(x^v(i,j)\right),\ V(i,j)=V\left(x(i,j)\right)$$

$$V^h(i+1,j)=V^h\left(x(i+1,j)\right),\ V^v(i,j+1)=V^v\left(x(i,j+1)\right)$$

Thus, it is easy to get that

$$V^h(i+1,j)+V^v(i,j+1)<\alpha\left(V^h(i,j)+V^v(i,j)\right) \tag{14}$$

Since for any nonnegative integer $D>z=\max(z_1,z_2)$, we have $V^h(0,D)=V^v(D,0)=0$, then summing up both sides of (14) from $D-1$ to $0$ with respect to $j$ and $0$ to $D-1$ with respect to $i$, for any nonnegative integer $D>D'\ge z=\max(z_1,z_2)$, one gets

$$\begin{aligned}\sum_{i+j=D}V(i,j)= &\ V^h(0,D)+V^h(1,D-1)+V^h(2,D-2)+\cdots+V^h(D-1,1)+V^h(D,0)\\ &+V^v(0,D)+V^v(1,D-1)+V^v(2,D-2)+\cdots+V^v(D-1,1)+V^v(D,0)\\ &<\alpha\left(V^h(0,D-1)+V^v(0,D-1)+V^h(1,D-2)+V^v(1,D-2)\right.\\ &\quad\left.+\cdots+V^h(D-1,0)+V^v(D-1,0)\right)\\ &=\alpha\sum_{i+j=D-1}V(i,j)<\ldots<\alpha^{D-D'}\sum_{i+j=D'}V(i,j)\end{aligned}$$

The proof is completed.

**Lemma 4** Consider the system (8), for a given positive constant $\beta>1$, if there exist positive definite symmetric matrices $P=diag\{P_h,P_v\}$ , $Q=diag\{Q_h,Q_v\}$ with appropriate dimensions, such that

$$\begin{bmatrix} Q-\beta P & 0 & A^T P \\ * & -\Lambda_2 Q & A_d^T P \\ * & * & -P \end{bmatrix}<0 \tag{15}$$

where $\Lambda_2=diag\left\{\beta^{d_h}I_h,\beta^{d_v}I_v\right\}$.

then, along the trajectory of systems (8), there holds the following inequality

$$\sum_{i+j=D}V(i,j)<\beta^{D-D'}\sum_{i+j=D'}V(i,j) \tag{16}$$

where $D' \geq z$.

**Proof** Following the similar proof of Lemma 3, Lemma 4 can be derived, it is omitted here.

**Remark 4** Lemmas 3 and 4 provide the methods for the estimation of Lyapunov functional candidate which will be used to design the controller for the 2D discrete switched system under asynchronous switching.

### 3.2. Controller design

Consider systems (1), under asynchronous switching controller $u(i,j) = K^{\sigma'(m)} x(i,j)$, then the corresponding closed-loop system is given by

$$\begin{bmatrix} x^h(i+1,j) \\ x^v(i,j+1) \end{bmatrix} = \left( \hat{A}^{\sigma(m)} + \hat{B}^{\sigma(m)} K^{\sigma'(m)} \right) \begin{bmatrix} x^h(i,j) \\ x^v(i,j) \end{bmatrix} + A_d^{\sigma(m)} \begin{bmatrix} x^h(i-d_h,j) \\ x^v(i,j-d_v) \end{bmatrix} \tag{17}$$

Suppose that the $k-th$ subsystem is activated at the switching instant $m_\kappa$, the $l-th$ subsystem is activated at the switching instant $m_{\kappa+1}$, the corresponding switching controller is activated at the switching instant $m_\kappa + \Delta m_\kappa$, $m_{\kappa+1} + \Delta m_{\kappa+1}$, respectively.

Let $T^+(z,D)$ denote the total mismatched period during $[z,D)$, $T^-(z,D)$ denote the total matched period during $[z,D)$, then we can get the following result.

**Theorem 1** Consider system (1), for given positive constants $\alpha < 1$, $\beta > 1$, if there exists positive definite symmetric matrices $X^k = diag\{X_h^k, X_v^k\}$, $Y^k = diag\{Y_h^k, Y_v^k\}$, $X^{kl} = diag\{X_h^{kl}, X_v^{kl}\}$, $Y^{kl} = diag\{Y_h^{kl}, Y_v^{kl}\}$ with appropriate dimensions, and positive scalars $\varepsilon_k$, $\varepsilon_{kl}$ such that, for $k,l \in \underline{N}, k \neq l$, the following inequalities (18) and (19) hold.

$$\begin{bmatrix} -\alpha X^k & 0 & \left(A^k X^k + B^k W^k\right)^T & X^k & \left(E_1^k X^k + E_3^k W^k\right)^T \\ * & -\Lambda_1 Y^k & \left(A_d^k Y^k\right)^T & 0 & \left(E_2^k Y^k\right)^T \\ * & * & -X^k + \varepsilon_k H^k H^{kT} & 0 & 0 \\ * & * & * & -Y^k & 0 \\ * & * & * & * & -\varepsilon_k I \end{bmatrix} < 0 \tag{18}$$

$$\begin{bmatrix} -\beta X^{kl} & 0 & \left(A^l X^{kl} + B^l K^k X^{kl}\right)^T & X^{kl} & \left(E_1^l X^{kl} + E_3^l K^k X^{kl}\right)^T \\ * & -\Lambda_2 Y^{kl} & \left(A_d^l Y^{kl}\right)^T & 0 & \left(E_2^l Y^{kl}\right)^T \\ * & * & -X^{kl} + \varepsilon_{kl} H^l H^{lT} & 0 & 0 \\ * & * & * & -Y^{kl} & 0 \\ * & * & * & * & -\varepsilon_{kl} I \end{bmatrix} < 0 \tag{19}$$

Then, under the switching controller

$$u(i, j) = K^{\sigma'(m)} x(i, j), \quad K^k = W^k X_k^{-1} \tag{20}$$

and the following average dwell time scheme

$$\frac{T^-(z, D)}{T^+(z, D)} \geq \frac{\lambda^+ + \lambda^*}{\lambda^- - \lambda^*}, \quad \tau_a > \tau_a^* = \frac{\ln \mu_1 \mu_2}{\lambda^*} \tag{21}$$

the corresponding closed-loop system is exponentially stable, where $\lambda^- = -\ln \alpha$, $\lambda^+ = \ln \beta$, $0 < \lambda^* < \lambda^-$, $\mu = \left(\frac{\alpha}{\beta}\right)^{\bar{d}}$, $\bar{d} = \max\{d_h, d_v\}$, $\mu_1 \mu_2 \mu \geq 1$ satisfying

$$X_l^{-1} \leq \mu_1 X_{kl}^{-1}, \quad X_{kl}^{-1} \leq \mu_2 X_k^{-1}, \quad Y_l^{-1} \leq \mu_1 Y_{kl}^{-1}, \quad Y_{kl}^{-1} \leq \mu_2 \mu Y_k^{-1} \tag{22}$$

**Proof** When $D \in [m_\kappa + \Delta m_\kappa, m_{\kappa+1})$, the closed-loop system (17) can be written as

$$\begin{bmatrix} x^h(i+1, j) \\ x^v(i, j+1) \end{bmatrix} = \left(\hat{A}^k + \hat{B}^k K^k\right) \begin{bmatrix} x^h(i, j) \\ x^v(i, j) \end{bmatrix} + \hat{A}_d^k \begin{bmatrix} x^h(i - d_h, j) \\ x^v(i, j - d_v) \end{bmatrix} \tag{23}$$

For the system, we consider the following Lyapunov function candidate

$$V_k(x(i, j)) = V_k^h(x^h(i, j)) + V_k^v(x^v(i, j)) \tag{24}$$

where

$$V_k^h(x^h(i, j)) = \sum_{g=1}^{2} V_{gk}^h(x^h(i, j))$$

$$V_{1k}^h(x^h(i, j)) = x^h(i, j)^T P_h^k x^h(i, j)$$

$$V_{2k}^h(x^h(i, j)) = \sum_{r=i-d_h}^{i-1} x^h(r, j)^T Q_h^k x^h(r, j) \alpha^{i-r-1}$$

$$V_k^v(x^v(i, j)) = \sum_{g=1}^{2} V_{gk}^v(x^v(i, j))$$

$$V_{1k}^{v}\left(x^{v}\left(i,j\right)\right)=x^{v}\left(i,j\right)^{T}P_{v}^{k}x^{v}\left(i,j\right)$$

$$V_{2k}^{v}\left(x^{v}\left(i,j\right)\right)=\sum_{t=j-d_{v}}^{j-1}x^{v}\left(i,s\right)^{T}Q_{v}^{k}x^{v}\left(i,s\right)\alpha^{j-t-1}$$

By Lemma 3, we know that if there exist positive definite symmetric matrices $P^{k}=diag\left\{P_{h}^{k},P_{v}^{k}\right\}$, $Q^{k}=diag\left\{Q_{h}^{k},Q_{v}^{k}\right\}$ with appropriate dimensions such that

$$\begin{bmatrix} Q^{k}-\alpha P^{k} & 0 & \left(\hat{A}^{k}+\hat{B}^{k}K^{k}\right)^{T}P^{k} \\ * & -\Lambda_{1}Q^{k} & \hat{A}_{d}^{kT}P^{k} \\ * & * & -P^{k} \end{bmatrix}<0 \tag{25}$$

then the following inequality holds

$$\sum_{i+j=D}V_{k}\left(i,j\right)<\alpha^{D-D^{k}}\sum_{i+j=D^{k}}V_{k}\left(i,j\right) \tag{26}$$

where $D^{k}=m_{\kappa}+\Delta m_{\kappa}$.

When $D\in\left[m_{\kappa+1},m_{\kappa+1}+\Delta m_{\kappa+1}\right)$, the closed-loop system (17) can be written as

$$\begin{bmatrix} x^{h}\left(i+1,j\right) \\ x^{v}\left(i,j+1\right) \end{bmatrix}=\left(\hat{A}^{l}+\hat{B}^{l}K^{k}\right)\begin{bmatrix} x^{h}\left(i,j\right) \\ x^{v}\left(i,j\right) \end{bmatrix}+\hat{A}_{d}^{l}\begin{bmatrix} x^{h}\left(i-d_{h},j\right) \\ x^{v}\left(i,j-d_{v}\right) \end{bmatrix}$$

We consider the following Lyapunov function candidate

$$V_{kl}\left(x\left(i,j\right)\right)=V_{kl}^{h}\left(x^{h}\left(i,j\right)\right)+V_{kl}^{v}\left(x^{v}\left(i,j\right)\right) \tag{27}$$

where

$$V_{kl}^{h}\left(x^{h}\left(i,j\right)\right)=\sum_{g=1}^{2}V_{gkl}^{h}\left(x^{h}\left(i,j\right)\right)$$

$$V_{1kl}^{h}\left(x^{h}\left(i,j\right)\right)=x^{h}\left(i,j\right)^{T}P_{h}^{kl}x^{h}\left(i,j\right)$$

$$V_{2kl}^{h}\left(x^{h}\left(i,j\right)\right)=\sum_{r=i-d_{h}}^{i-1}x^{h}\left(r,j\right)^{T}Q_{h}^{kl}x^{h}\left(r,j\right)\alpha^{i-r-1}$$

$$V_{kl}^{v}\left(x^{v}\left(i,j\right)\right)=\sum_{g=1}^{2}V_{gkl}^{v}\left(x^{v}\left(i,j\right)\right)$$

$$V_{1kl}^{v}\left(x^{v}\left(i,j\right)\right)=x^{v}\left(i,j\right)^{T}P_{v}^{kl}x^{v}\left(i,j\right)$$

$$V_{2kl}^{v}\left(x^{v}\left(i,j\right)\right)=\sum_{t=j-d_{v}}^{j-1}x^{v}\left(i,s\right)^{T}Q_{v}^{kl}x^{v}\left(i,s\right)\alpha^{j-t-1}$$

By Lemma 4, we know that if there exist positive definite symmetric matrices $P^{kl}=diag\left\{P_{h}^{kl},P_{v}^{kl}\right\}$, $Q^{kl}=diag\left\{Q_{h}^{kl},Q_{v}^{kl}\right\}$, with appropriate dimensions such that

$$\begin{bmatrix} Q^{kl}-\beta P^{kl} & 0 & \left(\hat{A}^{l}+\hat{B}^{l}K^{k}\right)^{T}P^{kl} \\ * & -\Lambda_{2}Q^{kl} & \hat{A}_{d}^{lT}P^{kl} \\ * & * & -P^{kl} \end{bmatrix}<0 \tag{28}$$

then the following inequality holds

$$\sum_{i+j=D}V_{kl}\left(i,j\right)<\beta^{D-D^{kl}}\sum_{i+j=D^{kl}}V_{kl}\left(i,j\right) \tag{29}$$

where $D^{kl}=m_{\kappa+1}$.

We consider the following piece-wise Lyapunov functional candidate for the closed-loop system (17)

$$V\left(i,j\right)=\begin{cases} x^{h}\left(i,j\right)^{T}P_{h}^{k}x^{h}\left(i,j\right)+x^{v}\left(i,j\right)^{T}P_{v}^{k}x^{v}\left(i,j\right) \\ +\sum_{r=i-d_{h}}^{i-1}x^{h}\left(r,j\right)^{T}Q_{h}^{k}x^{h}\left(r,j\right)\alpha^{i-r-1}+\sum_{t=j-d_{v}}^{j-1}x^{v}\left(i,s\right)^{T}Q_{v}^{k}x^{v}\left(i,s\right)\alpha^{j-t-1}, \\ \qquad D\in\left[m_{\pi}+\Delta m_{\pi},m_{\pi+1}\right),\quad \pi=1,2,...\kappa... \\ x^{h}\left(i,j\right)^{T}P_{h}^{kl}x^{h}\left(i,j\right)+x^{v}\left(i,j\right)^{T}P_{v}^{kl}x^{v}\left(i,j\right) \\ +\sum_{r=i-d_{h}}^{i-1}x^{h}\left(r,j\right)^{T}Q_{h}^{kl}x^{h}\left(r,j\right)\alpha^{i-r-1}+\sum_{t=j-d_{v}}^{j-1}x^{v}\left(i,s\right)^{T}Q_{v}^{kl}x^{v}\left(i,s\right)\alpha^{j-t-1}, \\ \qquad D\in\left[m_{\pi},m_{\pi}+\Delta m_{\pi}\right),\quad \pi=1,2,...\kappa... \end{cases} \tag{30}$$

Now Let $\chi=N_{\sigma(m)}(z,D)$ denote the switch number of $\sigma\left(m\right)$ on an interval $\left(z,D\right)$, and let $m_{\kappa-\chi+1}<\ m_{\kappa-\chi+2}<\cdots<m_{\kappa-1}<m_{\kappa}$ denote the switching points of $\sigma\left(m\right)$ over the interval $\left(z,D\right)$, then the switching points of $\sigma'\left(m\right)$ can be denoted as follows

$$m_{\kappa-\chi+1}+\Delta m_{\kappa-\chi+1}<m_{\kappa-\chi+2}+\Delta m_{\kappa-\chi+2}<\cdots<m_{\kappa-1}+\Delta m_{\kappa-1}<m_{\kappa}+\Delta m_{\kappa}$$

When $m_{\kappa-\upsilon+1}>z\geq m_{\kappa-\upsilon}+\Delta m_{\kappa-\upsilon}$, for $D>m_{\kappa}+\Delta m_{\kappa}$, using the condition (22) and (30), we

can get

$$\sum_{i+j=m_\kappa+\Delta m_\kappa} V(i,j) \le \mu_1 \sum_{i+j=(m_\kappa+\Delta m_\kappa)^-} V(i,j), \quad \sum_{i+j=m_\kappa} V(i,j) \le \mu_2 \sum_{i+j=(m_\kappa)^-} V(i,j)$$

where $(m_\kappa)^-$ and $(m_\kappa+\Delta m_\kappa)^-$ satisfies the following conditions

$$0 < m_\kappa - (m_\kappa)^- < \lambda, \quad 0 < m_\kappa + \Delta m_\kappa - (m_\kappa+\Delta m_\kappa)^- < \lambda,$$

where $\lambda$ is a sufficient small positive constant.

$$\begin{aligned}
&\sum_{i+j=D} V(i,j) < \alpha^{D-m_\kappa-\Delta m_\kappa} \sum_{i+j=m_\kappa+\Delta m_\kappa} V(i,j) \le \mu_1 \alpha^{D-m_\kappa-\Delta m_\kappa} \sum_{i+j=(m_\kappa+\Delta m_\kappa)^-} V(i,j) \\
&< \mu_1 \alpha^{D-m_\kappa-\Delta m_\kappa} \beta^{\Delta m_\kappa} \sum_{i+j=m_\kappa} V(i,j) \le \mu_1\mu_2 \alpha^{D-m_\kappa-\Delta m_\kappa} \beta^{\Delta m_\kappa} \sum_{i+j=(m_\kappa)^-} V(i,j) < \dots \\
&< (\mu_1\mu_2)^v \alpha^{D-m_\kappa-\Delta m_\kappa+m_\kappa-m_{\kappa-1}-\Delta m_{\kappa-1}+\dots+m_{\kappa-\upsilon+2}-m_{\kappa-\upsilon+1}-\Delta m_{\kappa-\upsilon+1}} \beta^{\Delta m_\kappa+\Delta m_{\kappa-1}+\dots+\Delta m_{\kappa-\upsilon+1}} \sum_{i+j=(m_{k-\upsilon+1})^-} V(i,j) \\
&< (\mu_1\mu_2)^v \alpha^{D-m_\kappa-\Delta m_\kappa+m_\kappa-m_{\kappa-1}-\Delta m_{\kappa-1}+\dots+m_{\kappa-\upsilon+2}-m_{\kappa-\upsilon+1}-\Delta m_{\kappa-\upsilon+1}+m_{\kappa-\upsilon+1}-z} \beta^{\Delta m_\kappa+\Delta m_{\kappa-1}+\dots+\Delta m_{\kappa-\upsilon+1}} \sum_{i+j=z} V(i,j) \\
&= (\mu_1\mu_2)^v \alpha^{T^-(z,D)} \beta^{T^+(z,D)} \sum_{i+j=z} V(i,j) = e^{v\ln(\mu_1\mu_2)+T^-(z,D)\ln\alpha+T^+(z,D)\ln\beta} \sum_{i+j=z} V(i,j) \qquad (31)
\end{aligned}$$

When $m_{\kappa-\upsilon} < z < m_{\kappa-\upsilon} + \Delta m_{\kappa-\upsilon}$, we can also get

$$\begin{aligned}
&\sum_{i+j=D} V(i,j) < \alpha^{D-m_\kappa-\Delta m_\kappa} \sum_{i+j=m_\kappa+\Delta m_\kappa} V(i,j) < \mu_1 \alpha^{D-m_\kappa-\Delta m_\kappa} \sum_{i+j=(m_\kappa+\Delta m_\kappa)^-} V(i,j) \\
&< \mu_1 \alpha^{D-m_\kappa-\Delta m_\kappa} \beta^{\Delta m_\kappa} \sum_{i+j=m_\kappa} V(i,j) < \mu_1\mu_2 \alpha^{D-m_\kappa-\Delta m_\kappa} \beta^{\Delta m_\kappa} \sum_{i+j=(m_\kappa)^-} V(i,j) < \dots \\
&< (\mu_1\mu_2)^v \alpha^{D-m_\kappa-\Delta m_\kappa+m_\kappa-m_{\kappa-1}-\Delta m_{\kappa-1}+\dots+m_{\kappa-\upsilon+2}-m_{\kappa-\upsilon+1}-\Delta m_{\kappa-\upsilon+1}} \beta^{\Delta m_\kappa+\Delta m_{\kappa-1}+\dots+\Delta m_{\kappa-\upsilon+1}} \sum_{i+j=(m_{k-\upsilon+1})^-} V(i,j) \\
&< (\mu_1\mu_2)^v \alpha^{D-m_\kappa-\Delta m_\kappa+m_\kappa-m_{\kappa-1}-\Delta m_{\kappa-1}+\dots+m_{\kappa-\upsilon+2}-m_{\kappa-\upsilon+1}-\Delta m_{\kappa-\upsilon+1}+m_{\kappa-\upsilon+1}-m_{\kappa-\upsilon}-\Delta m_{\kappa-\upsilon}} \beta^{\Delta m_\kappa+\Delta m_{\kappa-1}+\dots+\Delta m_{\kappa-\upsilon+1}+\Delta m_{\kappa-\upsilon}+m_{\kappa-\upsilon}-z} \sum_{i+j=z} V(i,j) \\
&= \mu_1 (\mu_1\mu_2)^v \alpha^{T^-(z,D)} \beta^{T^+(z,D)} \sum_{i+j=z} V(i,j) = \mu_1 e^{v\ln(\mu_1\mu_2)+T^-(z,D)\ln\alpha+T^+(z,D)\ln\beta} \sum_{i+j=z} V(i,j) \qquad (32)
\end{aligned}$$

From (31) and (32), we have

$$\sum_{i+j=D} V(i,j) < \max\{\mu_1, 1\} e^{v\ln(\mu_1\mu_2)+T^-(z,D)\ln\alpha+T^+(z,D)\ln\beta} \sum_{i+j=z} V(i,j) \qquad (33)$$

According to definition 2, we know

$$v = N_{\sigma(m)}(z,D) \le N_0 + \frac{D-z}{\tau_a} \qquad (34)$$

From condition (31), we can obtain

$$-T^-(z,D)\lambda^- + T^+(z,D)\lambda^+ \le -\lambda^*(D-z) \qquad (35)$$

Then for (33), we can obtain

$$\sum_{i+j=D} V(i,j) < \max\{\mu_1,1\}(\mu_1\mu_2)^{N_0} e^{\left(\frac{\ln(\mu_1\mu_2)}{\tau_a}-\lambda^*\right)(D-z)} \sum_{i+j=z} V(i,j) \tag{36}$$

It follows that

$$\sum_{i+j=D} \left\|x(i,j)\right\|^2 \le (\zeta_1/\zeta_2)\max\{\mu_1,1\}(\mu_1\mu_2)^{N_0} e^{\left(\frac{\ln(\mu_1\mu_2)}{\tau_a}-\lambda^*\right)(D-z)} \sum_{i+j=z} \left\|x(i,j)\right\|_C^2 \tag{37}$$

where

$$\zeta_1 = \max_{k,l\in \underline{N},k\ne l}\left\{\lambda_{\max}\left(P^k\right)+\max(d_h,d_v)\lambda_{\max}\left(Q^k\right),\lambda_{\max}\left(P^{kl}\right)+\max(d_h,d_v)\lambda_{\max}\left(Q^{kl}\right)\right\}$$

$$\zeta_2 = \min_{k,l\in \underline{N},k\ne l}\left\{\lambda_{\min}\left(P^k\right),\lambda_{\min}\left(P^{kl}\right)\right\}$$

$$\sum_{i+j=z}\left\|x(i,j)\right\|_C^2 \triangleq \sup_{\substack{-d_h<\theta_h\le 0,\\ -d_v<\theta_v\le 0}} \sum_{i+j=z}\left\{\left\|x^h(i-\theta_h,j)\right\|^2+\left\|x^v(i,j-\theta_v)\right\|^2\right\}$$

So when the condition (21) is satisfied, $\sum_{i+j=D}\left\|x(i,j)\right\|^2$ will tend to zero exponentially, which means the corresponding closed-loop 2D discrete switched system is exponentially stable.

Denoting $X^k=\left(P^k\right)^{-1}$, $Y^k=\left(Q^k\right)^{-1}$, then it is easy to get that $\left(X^k\right)^T=X^k$, $\left(Y^k\right)^T=Y^k$. Using $diag\left\{X^k,Y^k,X^k\right\}$ to pre- and post-multiply the left of (25), respectively, and denoting $W^k=K^kX^k$, then by lemma 1, it follows that (38) and (25) are equivalent.

$$\begin{bmatrix} -\alpha X^k & 0 & \left(\hat{A}^kX^k+\hat{B}^kW^k\right)^T & X^k \\ * & -\Lambda_1 Y^k & \left(\hat{A}_d^kY^k\right)^T & 0 \\ * & * & -X^k & 0 \\ * & * & * & -Y^k \end{bmatrix} < 0 \tag{38}$$

By substituting (2) into (38), we get the following inequality

$$T = T_0 + T_1 < 0 \tag{39}$$

where

$$\begin{bmatrix} -\alpha X^k & 0 & \left(A^k X^k + B^k W^k\right)^T & X^k \\ * & -\Lambda_1 Y^k & \left(A_d^k Y^k\right)^T & 0 \\ * & * & -X^k & 0 \\ * & * & * & -Y^k \end{bmatrix} < 0$$

$$T_1 = \begin{bmatrix} \left(E_1^k X^k + E_3^k W^k\right)^T \\ \left(E_2^k Y^k\right)^T \\ 0 \\ 0 \end{bmatrix} F^T \begin{bmatrix} 0 \\ 0 \\ H^k \\ 0 \end{bmatrix}^T + \begin{bmatrix} 0 \\ 0 \\ H^k \\ 0 \end{bmatrix} F \begin{bmatrix} \left(E_1^k X^k + E_3^k W^k\right)^T \\ \left(E_2^k Y^k\right)^T \\ 0 \\ 0 \end{bmatrix}^T$$

By lemma 2, we can get

$$T_1 \le \varepsilon_k^{-1} \begin{bmatrix} \left(E_1^k X^k + E_3^k W^k\right)^T \\ \left(E_2^k Y^k\right)^T \\ 0 \\ 0 \end{bmatrix} \begin{bmatrix} \left(E_1^k X^k + E_3^k W^k\right)^T \\ \left(E_2^k Y^k\right)^T \\ 0 \\ 0 \end{bmatrix}^T + \varepsilon_k \begin{bmatrix} 0 \\ 0 \\ H^k \\ 0 \end{bmatrix} \begin{bmatrix} 0 \\ 0 \\ H^k \\ 0 \end{bmatrix}^T \tag{40}$$

Applying lemma 1 again, we obtain that (39) is equivalent to (18).

Similarly, denoting $X^{kl} = \left(P^{kl}\right)^{-1}$, $Y^{kl} = \left(Q^{kl}\right)^{-1}$, and using $diag\left\{X^{kl}, Y^{kl}, X^{kl}\right\}$ to pre- and post-multiply the left of (28), respectively, then substituting (2) into (28) and applying lemma 1 and 2, it is easy to get that (28) is equivalent to (19). The proof is completed.

**Remark 4** In Theorem 1, we propose sufficient conditions for the existence of a state feedback controller such that the considered 2D discrete switched system (1) under asynchronous switching is exponentially stable. It is worth noting that this condition is obtained by using the average dwell-time approach. Here, $\alpha$ plays a key role in controlling the rate of decaying of the closed-loop system in the matched period, and $\beta$ plays a key role in controlling the rate of increasing of the closed-loop system in the mismatched period, respectively.

**Remark 5** We would like to point out that there is still enough room to improve the result. Because time delays always exist in a time-varying fashion in real-time systems, so we can also consider 2D discrete switched systems with time-varying delays. The results can be improved by combining with free-weighting matrices to reduce conservatism [10], which is shown to lead to much less conservative results than most existing literature.

**Remark 6** It is noticed that (18) and (19) are mutually dependent. Therefore, we can firstly solve the LMI (18) to obtain the solutions of matrices $X^k$, $Y^k$ and $W^k$. Then (19) can be transformed into the LMI by substituting $K^k = W^k X_k^{-1}$ into (19). By adjusting the parameters $\alpha$, $\beta$, we can find the feasible solutions of $X^k$, $Y^k$, $W^k$, $X^{kl}$ and $Y^{kl}$ such that (18) and (19) hold.

**Remark 7** When $\Delta m_\kappa = 0, \kappa=1,2,\ldots$, Theorem 1can be reduced to the result of system (1) under synchronous switching.

When $A_d^{\sigma(m)} = 0$, the system (17) degenerates to the following system

$$\begin{bmatrix} x^h(i+1,j) \\ x^v(i,j+1) \end{bmatrix} = \left(\hat{A}^{\sigma(m)} + \hat{B}^{\sigma(m)} K^{\sigma'(m)}\right) \begin{bmatrix} x^h(i,j) \\ x^v(i,j) \end{bmatrix} \tag{41}$$

Then we can get the following corollary from Theorem 1.

**Corollary 1** Consider system (41), for given positive constants $\alpha < 1$, $\beta > 1$, if there exists positive definite symmetric matrices $X^k = diag\left\{X_h^k, X_v^k\right\}$, $X^{kl} = diag\left\{X_h^{kl}, X_v^{kl}\right\}$ ,with appropriate dimensions, and positive scalars $\varepsilon_k$, $\varepsilon_{kl}$ such that, for $k,l \in \underline{N}, k \neq l$, the following inequalities (42) and (43) hold.

$$\begin{bmatrix} -\alpha X^k & \left(A^k X^k + B^k W^k\right)^T & \left(E_1^k X^k + E_3^k W^k\right)^T \\ * & -X^k + \varepsilon_k H^k H^{kT} & 0 \\ * & * & -\varepsilon_k I \end{bmatrix} < 0 \tag{42}$$

$$\begin{bmatrix} -\beta X^{kl} & \left(A^l X^{kl} + B^l K^k X^{kl}\right)^T & \left(E_1^l X^{kl} + E_3^l K^k X^{kl}\right)^T \\ * & -X^{kl} + \varepsilon_{kl} H^l H^{lT} & 0 \\ * & * & -\varepsilon_{kl} I \end{bmatrix} < 0 \tag{43}$$

Then, under the switching controller

$$u(i,j) = K^{\sigma'(m)} x(i,j), \quad K^k = W^k X_k^{-1} \tag{44}$$

and the following average dwell time scheme

$$\frac{T^-(z,D)}{T^+(z,D)} \geq \frac{\lambda^+ + \lambda^*}{\lambda^- - \lambda^*}, \quad \tau_a > \tau_a^* = \frac{\ln \mu_1 \mu_2}{\lambda^*} \tag{45}$$

the corresponding closed-loop system is exponentially stable, where

$\lambda^- = -\ln\alpha$, $\lambda^+ = \ln\beta$, $0 < \lambda^* < \lambda^-$, $\mu_1\mu_2 \geq 1$ satisfying

$$X_l^{-1} \leq \mu_1 X_{kl}^{-1}, \quad X_{kl}^{-1} \leq \mu_2 X_k^{-1} \tag{46}$$

The procedure of the controller design for system (1) can be given as follows:

**Step 1** Given system matrices and positive constant $0 < \alpha < 1$, by solving the LMI (18), we can get the feasible solution of matrices $X^k$, $Y^k$ and $W^k$, and positive scalar $\varepsilon_k$, then the controller can be obtained by $K^k = W^k X_k^{-1}$.

**Step 2** Substitute $K^k$ into (19), then by adjusting the parameter $\alpha$, $\beta$, we can find the feasible solutions $X^k$, $Y^k$, $X^{kl}$, $Y^{kl}$ , $W^k$, $\varepsilon_k$ and $\varepsilon_{kl}$ such that (19) holds.

**Step 3** To minish the infimum of the average dwell time $\tau_a^*$, we can obtain the minimum of $\mu_1$ and $\mu_2$ satisfying (22) and $\mu_1\mu_2\mu \geq 1$.

**Step 4** Take $\lambda^*$ satisfying $0 < \lambda^* < \lambda^-$, and compute the value of $T^-(z,D)/T^+(z,D)$ and $\tau_a^*$ by (21).

## 4. Numerical example

In this section, we present an example to illustrate the effectiveness of the proposed approach.

$$A_1=\begin{bmatrix}1 & 1.5\\ 1 & 0.5\end{bmatrix},\ A_{d1}=\begin{bmatrix}-0.15 & 0\\ -0.1 & -0.12\end{bmatrix},\ B_1=\begin{bmatrix}-4.5 & 0\\ 1 & -3\end{bmatrix},\ H_1=\begin{bmatrix}0.2 & 0.15\\ 0.1 & 0.2\end{bmatrix},$$

$$E_1^1=\begin{bmatrix}0.1 & 0.15\\ 0.1 & 0\end{bmatrix},\ E_2^1=\begin{bmatrix}0.1 & 0\\ 0.1 & 0.2\end{bmatrix},\ F_1=\begin{bmatrix}\sin\left(0.5\pi\left(i+j\right)\right) & 0\\ 0 & \sin\left(0.5\pi\left(i+j\right)\right)\end{bmatrix},$$

$$E_3^1=\begin{bmatrix}0.15 & 0\\ 0.13 & 0.12\end{bmatrix},\ A_2=\begin{bmatrix}1 & 2\\ 1 & 1\end{bmatrix},\ A_{d2}=\begin{bmatrix}-0.1 & 0.2\\ 0 & -0.2\end{bmatrix},\ B_2=\begin{bmatrix}-5 & 1\\ -1 & -3\end{bmatrix},$$

$$H_2=\begin{bmatrix}0.2 & 0.25\\ 0.2 & 0.3\end{bmatrix},\ E_1^2=\begin{bmatrix}0.1 & 0.2\\ 0.2 & 0.1\end{bmatrix},\ E_2^2=\begin{bmatrix}0.2 & 0.1\\ 0.2 & 0.1\end{bmatrix},\ E_3^2=\begin{bmatrix}0.12 & 0.15\\ 0.12 & 0.1\end{bmatrix},$$

$$F_2=\begin{bmatrix}\cos\left(0.5\pi\left(i+j\right)\right) & 0\\ 0 & \cos\left(0.5\pi\left(i+j\right)\right)\end{bmatrix},\ d_h=2,\ d_v=3.$$

The bounded conditions are as follows

$$x^h\left(0,j\right)=\begin{cases}10, & 0\le j\le 20\\ 0, & j>20\end{cases},\quad x^v\left(i,0\right)=\begin{cases}6, & 0\le i\le 20\\ 0, & i>20\end{cases}.$$

where state dimensions $n_1=1$, $n_2=1$.

Take $\alpha=0.6$, $\beta=1.2$, then solving the matrix inequalities (18) in theorem 1, gives rise to

$$X^1=\begin{bmatrix}0.4300 & 0.0080\\ 0.0080 & 0.4365\end{bmatrix},\ X^2=\begin{bmatrix}0.4180 & -0.0451\\ -0.0451 & 0.3841\end{bmatrix},$$

$$Y^1=\begin{bmatrix}1.0619 & -0.0738\\ -0.0738 & 1.0681\end{bmatrix},\ Y^2=\begin{bmatrix}1.0295 & -0.0395\\ -0.0395 & 0.8684\end{bmatrix},$$

$$W^1=\begin{bmatrix}0.0991 & 0.1475\\ 0.1789 & 0.1255\end{bmatrix},\ W^2=\begin{bmatrix}0.0840 & 0.1560\\ 0.0944 & 0.0590\end{bmatrix}$$

$$\varepsilon_1=0.7882,\ \varepsilon_2=0.9010$$

Then $K^1$, $K^2$ can be obtained by $K^k=W^kX_k^{-1}$

$$K^1=\begin{bmatrix}0.2243 & 0.3338\\ 0.4107 & 0.2800\end{bmatrix},\ K^2=\begin{bmatrix}0.2481 & 0.4354\\ 0.2454 & 0.1823\end{bmatrix}.$$

Substituting $K^1$ and $K^2$ into (19), and solving the matrix inequalities (19), we can get the following matrices

$$X^{12}=\begin{bmatrix}0.6059 & -0.1005\\ -0.1005 & 0.5325\end{bmatrix},\quad X^{21}=\begin{bmatrix}0.5805 & -0.0926\\ -0.0926 & 0.5677\end{bmatrix},$$

$$Y^{12}=\begin{bmatrix}1.0435 & -0.0451\\ -0.0451 & 0.9365\end{bmatrix},\quad Y^{21}=\begin{bmatrix}1.0422 & -0.0447\\ -0.0447 & 0.9670\end{bmatrix}.$$

$$\varepsilon_{12}=0.9987,\quad \varepsilon_{21}=1.0404$$

Then, we can get the minimum $\mu_1=1.5694$, $\mu_2=9.1561$, from (22). Let $\frac{T^-(z,D)}{T^+(z,D)}\geq\frac{\lambda^++\lambda^*}{\lambda^--\lambda^*}=8$, it is easy to get that $\lambda^*=0.4338$. Then from (21), it follows that $\tau_a^*=6.15$. Choosing $\tau_a=6.5$, and taking $\Delta m_\kappa=2, \kappa=1,2,\ldots$, the trajectories of the states $x^h(i,j)$ and $x^v(i,j)$ are shown in Fig. 1 and Fig. 2, respectively, and the system switching signal $\sigma(m)$ and the controller switching signal $\sigma'(m)$ are shown in Fig. 3. One can notice that the states of the closed-loop system converge to zero under the asynchronous switching.

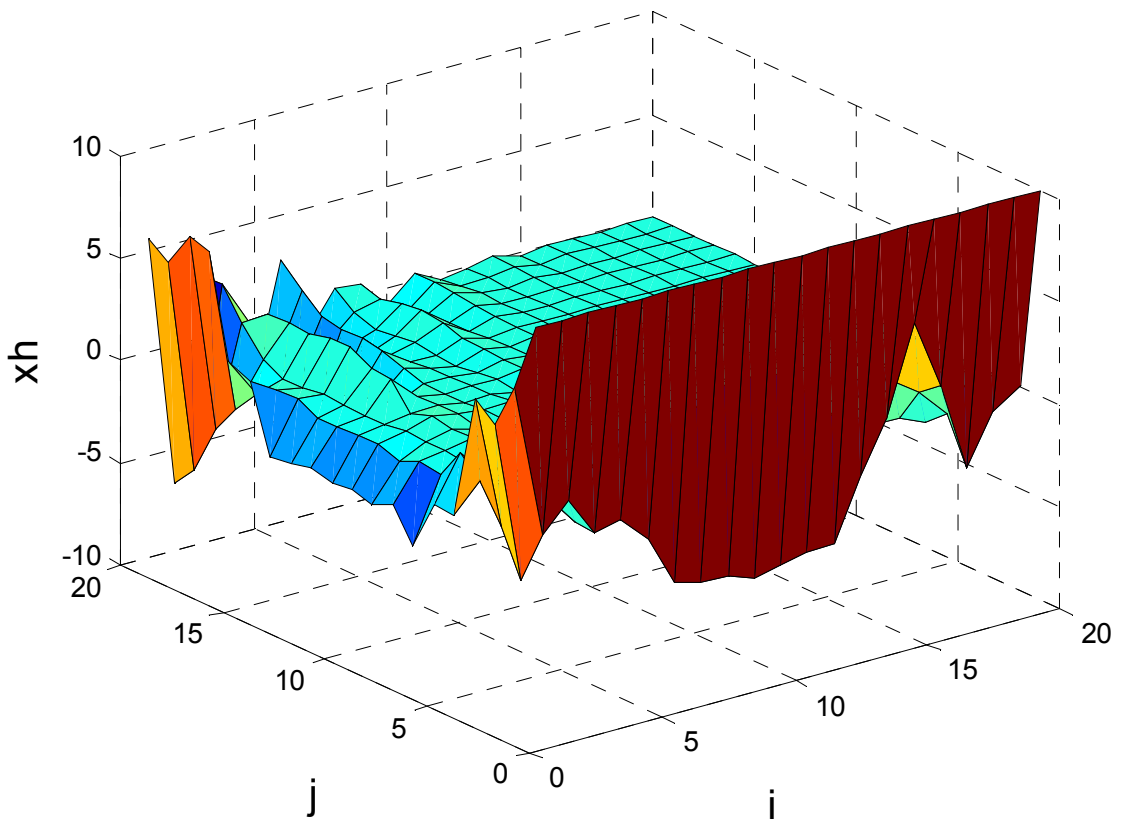

Fig. 1 The trajectory of the state $x^h(i,j)$

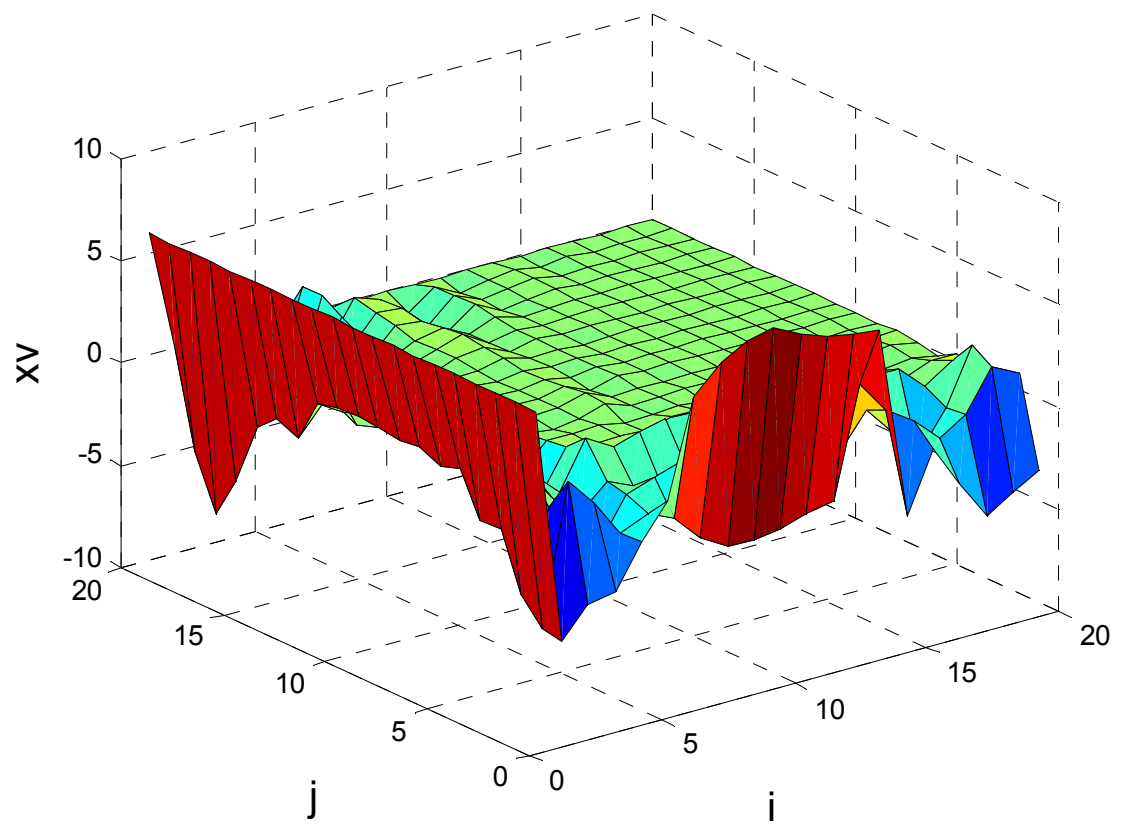

Fig. 2 The trajectory of the state $x^v(i,j)$

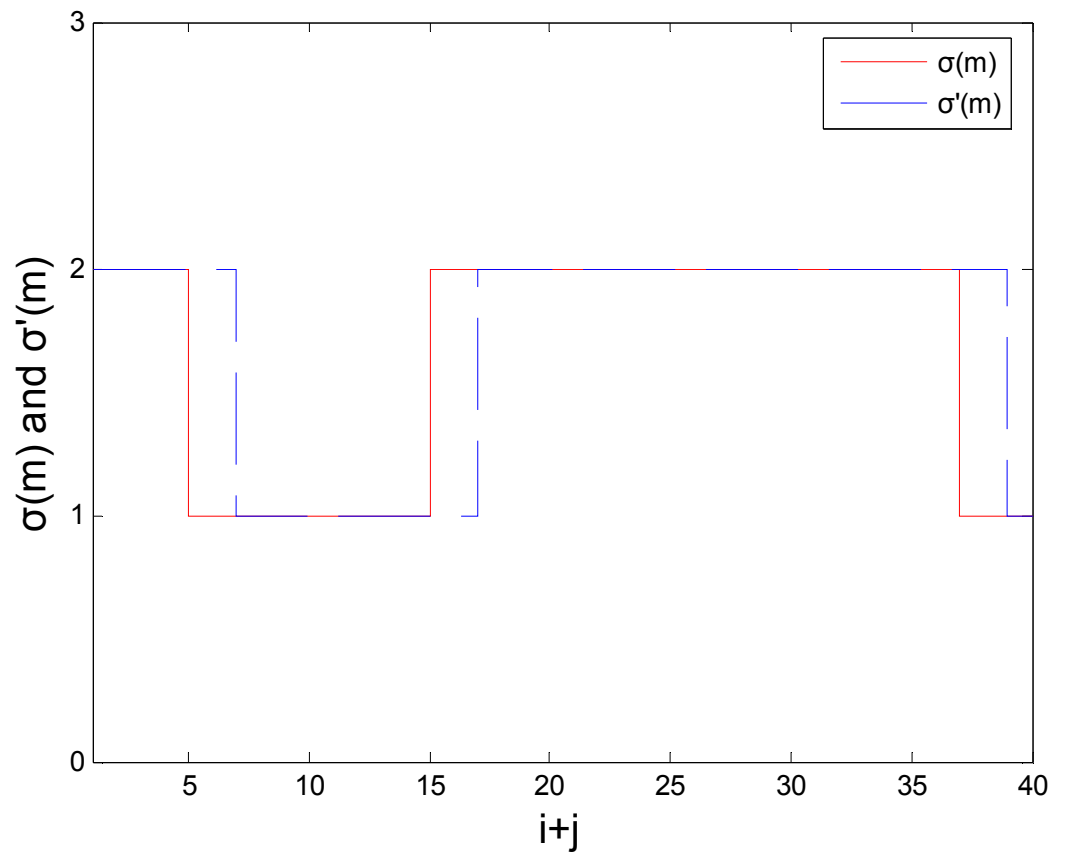

Fig. 3 Switching signal

# 5. Conclusion

This paper has investigated the problems of stability analysis and stabilization for a class of 2D discrete switched systems with constant state delays under asynchronous switching. A state feedback controller is proposed to stabilize this system, and the dwell time approach is utilized for the stability analysis. Sufficient conditions for the existence of such controller are formulated in terms of a set of LMI. An illustrative example is also given to illustrate the applicability of the proposed approach. Our future work will focus on extending the proposed design method to other problem such as robust $H_\infty$ control for 2D discrete switched systems with time-varying state delays under asynchronous switching.